\documentclass[a4paper,11pt]{amsart}
\usepackage{amssymb}
\usepackage[english]{babel}
\usepackage[latin1]{inputenc}
\usepackage{graphicx}

\theoremstyle{plain}
\newcommand{\ep}{\epsilon}
\newcommand{\R}{\mathbb R}
\newcommand{\A}{\mathcal A}
\newcommand{\B}{\mathcal B}
\newcommand{\M}{\mathcal M}

\begin{document}
\title[Convexity preserving models]{Convexity preserving jump diffusion models for
option pricing}
\author[Erik Ekström and Johan Tysk]{Erik Ekström$^1$ and Johan Tysk$^{2,3}$}
\date{January 4, 2006}
\email{ekstrom@maths.manchester.ac.uk, Johan.Tysk@math.uu.se}
\thanks{$^1$ School of Mathematics, University of Manchester,
Sackville Street, Manchester M60 1QD, UK}
\thanks{$^2$ Department of Mathematics, Uppsala University,
Box 480, SE-75106 Uppsala, Sweden}
\thanks{$^3$ Partially supported by the Swedish Research Council (VR)}

\newtheorem{theorem}{Theorem}[section]
\newtheorem{lemma}[theorem]{Lemma}
\newtheorem{corollary}[theorem]{Corollary}
\newtheorem{proposition}[theorem]{Proposition}
\newtheorem{definition}[theorem]{Definition}
\newtheorem{hypothesis}[theorem]{Hypothesis}

\newenvironment{example}[1][Example]{\begin{trivlist}
\item[\hskip \labelsep {\bf Example}]}{\end{trivlist}}
\newenvironment{remark}[1][Remark]{\begin{trivlist}
\item[\hskip \labelsep {\bf Remark}]}{\end{trivlist}}

\begin{abstract}
We investigate which jump-diffusion models are convexity preserving.
The study of convexity preserving models is motivated by
monotonicity results for such models in the volatility and in the
jump parameters. We give a necessary condition for convexity to be
preserved in several-dimensional jump-diffusion models. This
necessary condition is then used to show that, within a large class
of possible models, the only convexity preserving models are the
ones with linear coefficients.
\end{abstract}

\maketitle

\section{Introduction}
A model for a set of stock prices is said to be convexity preserving if the
price of any convex European claim is convex as a function of the
underlying stock prices at all times prior to maturity. As is
well-known, this property is intimately connected to certain
monotonicity properties of the option price with respect to
volatility and other parameters of the model. Generally speaking, if
the option price is convex at all fixed times, then it is also
increasing in the volatility. This robustness property motivates the
study of convexity preserving models in finance.

Although these issues have been studied quite extensively during the
last decade, compare \cite{BGW}, \cite{EKJPS}, \cite{H}, \cite{JT1}
and \cite{K} for the case of one-dimensional diffusion models,
\cite{BR}, \cite{EJT} and \cite{JT2} for several-dimensional
diffusion models and \cite{ET} for one-dimensional jump-diffusion
models, the general picture for more advanced models is not yet
fully understood. In \cite{ET}, a sufficient condition for the
preservation of convexity in one-dimensional models with jumps is
provided. That condition, however, is not a necessary condition for
preservation of convexity. The main contribution of the present
paper is to give a necessary condition for convexity to be preserved
in jump-diffusion models in arbitrary dimensions. We also use this
necessary condition to show that, within a large class of possible
models, the only higher-dimensional convexity preserving models are
the ones with linear coefficients.

To analyze the convexity of an option price we employ the
characterization of the price as the unique viscosity solution to a
parabolic integro-differential equation
\begin{equation}\label{eq}u_t=\A u+\B u\end{equation}
with an appropriate terminal condition. In this equation, $\A$ is an
elliptic differential operator associated with the continuous
fluctuations of the stock price processes, whereas $\B$ is an
integro-differential operator associated with the jumps of the stock
price processes, compare Section~\ref{regularity} below.
Preservation of convexity of the solution to the equation (\ref{eq})
is dealt with using the notion of \emph{locally convexity
preserving} (LCP) operators. This concept was introduced and
analyzed in \cite{JT2}, and also used in \cite{EJT} and \cite{ET}.
Following these references, we show that the condition that
$\M=\A+\B$ is LCP at all points is necessary for convexity to be
preserved. We also show that $\M$ is LCP if and only if both $\A$
and $\B$ are LCP, i.e.
\[\left\{\begin{array}{c}
\mbox{The model is}\\
\mbox{convexity}\\ \mbox{preserving}\end{array}\right\} \implies
\left\{\begin{array} {c}\mbox{$\M$ is LCP}\\ \mbox{at all
points}\end{array}\right\} \iff \left\{\begin{array}{c}
\mbox{Both $\A$ and $\B$}\\ \mbox{are LCP}\\
\mbox{at all points}\end{array}\right\},\] compare
Theorem~\ref{hekscher}. Thus the characterization of LCP models
breaks down into two easier problems: (i) to describe which
diffusion models are LCP, and (ii) to describe which jump structures
are LCP. Issue (i) has been dealt with in \cite{JT2} and \cite{EJT},
and (ii) is dealt with in Theorem~\ref{wittbrattstrom} below.

The present paper is organized as follows. In Section~\ref{model} we
introduce the model and we motivate the study of convexity
preserving models by means of a monotonicity result. In
Section~\ref{regularity} we prove a technical regularity result
which is used in the sequel. In Section~\ref{LCP condition} we
introduce the LCP-condition, and we show that a model is convexity
preserving only if both the differential operator $\A$ and the
integro-differential operator $\B$ are LCP at all points. In
Section~\ref{LCP} we investigate which jump structures are LCP. This
investigation is continued in Section~\ref{finite} for models with
only a finite number of possible jump sizes, where we show that,
within a large class of possible models, all convexity preserving
models have linear diffusion coefficients and jump structures.

\section{The model and a monotonicity result}
\label{model}

We consider a market consisting of $n$ different stocks, the prices
of which are modeled by an $n$-dimensional stochastic process $X(t)$.
To specify $X$, let $W$ be an $n$-dimensional Brownian motion, and
let $v$ be a Poisson random measure on $[0,T]\times [0,1]$ with
intensity measure
\[q(dt,dz)=\lambda(t)\,dtdz,\]
where $\lambda$ is a deterministic function. Let $X$ be a
jump-diffusion satisfying the stochastic differential equation
\[dX=\beta(X(t-),t)\,dW+\int_0^1\phi(X(t-),t,z)\,\tilde v(dt,dz).\]
Here $\beta=(\beta_{ij})_{ij=1}^n$ is an $n\times n$-matrix,
$\phi=(\phi_1,...,\phi_n)$ is an $n$-dimensional vector and $\tilde
v$ is the compensated jump martingale measure defined by
\[\tilde v(dt,dz)=(v-q)(dt,dz).\]

\begin{remark}
In this model, jumps occur according to a Poisson process
\[Y(t):=\int_0^t\int_0^1 v(dt,dz)\]
with intensity $\lambda(t)$.
Associated with each jump is a label $z\in[0,1]$. The interpretation
is that a jump of $Y$ at time $t$ with label $z$ results in a jump
of size $\phi_i(X(t-),t,z)$ in the $i$th coordinate of $X$. Between
jumps, $X$ follows a continuous diffusion governed by the diffusion
coefficient $\beta(X(t),t)$ and the drift
$-\lambda\int_0^1\phi(X(t),t,z)\,dz$.
\end{remark}

We denote by $\R_+^n$ the space $(0,\infty)^n$, and we say that a
function $g:\R_+^n\to\R$ is of at most polynomial growth if there
exist constants $m$ and $C$ such that
\[\vert g(x)\vert\leq C\left(1+\vert x\vert^m\right)\]
for all $x\in\R_+^n$. Given a continuous contract function
$g:\R_+^n\to\R$ of at most polynomial growth, the price at time
$t\in[0,T_0]$ of an option paying $g(X(T_0))$ at time $T_0\in[0,
T]$, is $u(X(t),t)$. Here the function $u:\R_+^n\times[0,T_0]\to\R$
is given by
\[u(x,t):=E_{x,t}g(X(T_0)),\]
where the indices indicate that $X(t)=x$. Note that the conditions
(A1)-(A7) specified below and the polynomial growth of $g$ implies
that all moments of $X(T)$ are finite, compare Section 7 in
\cite{GS}. Consequently, the option value $u$ is finite.

\begin{remark}
We do not address the issue of how to choose an appropriate pricing
measure, but we rather assume that the model is specified directly
under the measure used for pricing options. Also note that there is
no discounting factor in the definition of the option price. Thus we
implicitly assume, without loss of generality for our purposes, that
all prices are quoted in terms of some bond price.
\end{remark}

We will throughout this paper work under the regularity
and growth assumptions (A1)-(A7). When specifying these, $D$ is a
positive constant and the H\"older exponent $\alpha$ is a constant
between 0 and 1.

\begin{itemize}
\item[(A1)]
For all $i,j=1,...,n,$ $\beta_{ij}:\R^+\times[0,T]\to\R$ is in
$C^{2,0}_\alpha\left(\R_+^n\times[0,T]\right)$.
\item[(A2)]
$\lambda\in C_\alpha([0,T])$.
\item[(A3)] $\sum_{j=1}^n\vert\beta_{ij}(x,t)\vert^2+
\vert\phi_i(x,t,z)\vert^2\leq Dx_i^2$ for all $i=1,...,n$.
\item[(A4)]
$\vert\beta(x,t)-\beta(y,t)\vert+\vert\phi(x,t,z)-\phi(y,t,z)\vert\leq
D\vert x-y\vert$.
\item[(A5)]
The matrix $\beta(x,t)$ is non-singular for all
$(x,t)\in\R_+^n\times [0,T]$.
\item[(A6)]
$\phi:\R_+^n\times[0,T]\times[0,1]\to\R^n$ is measurable, and
$\phi_i(\cdot,\cdot,z)\in C^{2,0}_\alpha(\R_+^n\times[0,T])$ with
the H\"older norms being uniformly bounded in the $z$-variable.
Moreover, for all $x$ and $t$ we have $\phi(x,t,z)\not=0$ for almost
all $z$.
\item[(A7)]
There exists $\gamma>-1$ such that $\phi_i(x,t,z)>\gamma x_i$ for
$i=1,...,n$.
\end{itemize}

\begin{definition}
A model $(\beta,\lambda,\phi)$ is convexity preserving on $[0,T]$ if
for all $t$ and $T_0$ with $0\leq t\leq T_0 \leq T$, the price
$u(x,t)$ of the option with pay-off $g(X(T_0))$ at $T_0$ is convex
in $x$ for any convex contract function $g$ of at most polynomial
growth.
\end{definition}

The main reason for studying preservation of convexity is, as
mentioned in the introduction, that convexity implies certain
monotonicity properties of the option price with respect to the
parameters of the model. The following result can be proven in a
similar way as Theorem~5.1 in \cite{ET}, in which the one-dimensional
case is treated.

\begin{theorem}
\label{mon} Let two models be given with parameters $(\beta,
\lambda, \phi)$ and $(\tilde\beta, \tilde\lambda, \tilde\phi)$,
respectively. Assume that
\begin{itemize}
\item[(i)]
$\tilde\lambda(t)\leq \lambda (t)$ for all $t\in[0,T]$;
\item[(ii)]
for each fixed $(x,t)\in\R_+^n\times[0,T]$ we have
$\tilde\phi(x,z,t)=k(z)\phi(x,t,z)$ for some $k(z)\in[0,1]$;
\item[(iii)]
for all $(x,t)\in\R_+^n\times[0,T]$ we have $\tilde\beta\tilde\beta^*\leq
\beta\beta^*$ as quadratic forms (here $\beta^*$ denotes the transpose of
$\beta$).
\end{itemize}
Also assume that at least one of the two models is convexity
preserving. Then, for any convex contract function $g$ of at most polynomial
growth we have
\[\tilde u(x,t)\leq u(x,t)\]
for all $(x,t)\in\R_+^n\times[0,T_0]$, where $\tilde u$ and $u$ are the two
option prices corresponding to the two different models.
\end{theorem}

\begin{remark}
Note that the most important special case of (ii) is when for all $x$, $t$
and $z$ there exists an $i\in\{1,...,n\}$ such that
$\phi_j(x,t,z)=\tilde\phi_j(x,t,z)=0$ for all $j\not= i$
(i.e. the case when at most one component of $X$ and
one component of $\tilde X$ jump at each given time), and
$\phi_i(x,t,z)/\tilde\phi_i(x,z,t)\geq 1$ if
$\tilde\phi_i(x,z,t)\not= 0$. Also note that the
condition (iii) is the same as the one used for diffusion models in
higher dimensions, compare \cite{EJT}.
\end{remark}

Since all one-dimensional diffusion models and all geometric
Brownian motions (not necessarily one-dimensional) are known to be
convexity preserving, see \cite{BGW}, \cite{EKJPS}, \cite{H} or
\cite{JT1}, the following consequence of Theorem~\ref{mon} is
immediate. It is the higher-dimensional analogue to a result in
\cite{BJP}.

\begin{corollary}
Assume that a model $(\beta, \lambda, \phi)$ and a convex contract
function $g$ of at most polynomial growth are given. If $n\geq 2$,
also assume that $\beta$ is the (possibly time-dependent) diffusion
matrix of a geometric Brownian motion, i.e. $\beta_{ij}
(x,t)=\gamma_{ij}(t)x_i$ for some deterministic functions
$\gamma_{ij}$. Then a lower bound for the corresponding option price
is given by the option price in the model $(\beta, 0, \phi)$ with no
jumps.
\end{corollary}

\section{Regularity of the value function}
\label{regularity}

Under weak conditions, see for example \cite{P}, the pricing
function $u$ is the unique viscosity solution of a parabolic
integro-differential equation
\begin{equation}
\label{pide} u_t+\mathcal Mu=0
\end{equation}
with terminal condition
\[u(x,T_0)=g(x).\]
In this equation, the operator $\mathcal M=\mathcal A+\mathcal B$,
where the second order differential operator $\mathcal A$ and the
integro-differential operator $\mathcal B$ are given by
\[\mathcal Au(x,t):=\sum_{i,j=1}^na_{ij}(x,t)u_{x_ix_j}(x,t)\]
and \begin{equation} \label{B} \mathcal
Bu(x,t):=\lambda(t)\int_0^1\big(u(x+\phi(x,t,z),t)-u(x,t)-\phi(x,t,z)\cdot\nabla
u(x,t)\big)\,dz
\end{equation}
respectively, and $a_{ij}$ are the coefficients of the matrix
$\beta\beta^*/2$. Under the assumptions (A1)-(A7), $u$ is not merely
a viscosity solution to (\ref{pide}), but it is also a classical
solution. Indeed, we show below that the value function $u$
is regular. The proof has certain similarities to the proof of
Theorem~3.2 in \cite{ET} in which the one-dimensional case is
treated. However, for the convenience of the reader, and
since the proofs differ at some points, we include the higher-dimensional
version in its full detail.

\begin{theorem}
\label{reg} Assume that $g\in C^4_\alpha(\R_+^n)$ and that $g$ is
globally Lipschitz continuous. Then $u\in C^{4,1}_\alpha
(\R_+^n\times[0,T_0])$.
\end{theorem}

\begin{proof}
Let $\psi:\R_+\to\R$ be a smooth function with $\psi^\prime>0$ such
that
\[\psi(s)=\left\{\begin{array}{cl}
s & \mbox{if }s>2\\
-1/s &\mbox{if }s<1.\end{array}\right.\] It follows from
It$\hat{\mbox{o}}$'s lemma that the $n$-dimensional stochastic process $Y(t)$,
where
\[Y_i(t)=\psi(X_i(t)),\] satisfies
\[dY_i=\tilde
b_i(Y(t-),t)\,dt+\sum_{j=1}^n\tilde\beta_{ij}(Y(t-),t)\,dW_j +
\int_0^1\tilde\phi_i(Y(t-),t,z)\,\tilde v(dt,dz)\] on $\R^n\times
[0,T]$. Here \begin{eqnarray*} \tilde b_i (y,t)&=&
\frac{1}{2}\psi^{\prime\prime}(\psi^{-1}(y_i))
\sum_{j=1}^n\beta_{ij}^2 (\psi^{-1}(y),t)\\
&&+\lambda(t)\int_0^1\left(\tilde\phi_i
(y,t,z)-\psi^\prime(\psi^{-1}(y_i))\phi_i(\psi^{-1}(y),t,z)\right)\,dz,
\end{eqnarray*}
\[\tilde\beta_{ij}(y,t)=\psi^\prime(\psi^{-1}(y_i))\beta_{ij}(\psi^{-1}(y),t),\]
and \[\tilde\phi_i(y,t,z)=\psi\left(\psi^{-1}(y_i)+
\phi_i(\psi^{-1}(y),t,z)\right)-y_i ,\] where
$\psi^{-1}(y):=(\psi^{-1}(y_1),...,\psi^{-1}(y_n))$. Now it is
straightforward to check that $\tilde b$, $\tilde \beta$ and $\tilde
\phi$ together with the initial condition $\tilde
g(y):=g(\psi^{-1}(y))$ satisfy the conditions (2.2)-(2.5) in
\cite{P}. According to Theorem~3.1 in \cite{P}, the function
$v(y,t):=u(\psi^{-1}(y),t)$ is a viscosity solution of the
integro-differential equation
\begin{equation}\label{vis}
\left\{\begin{array}{ll} v_t+\sum_{i,j=1}^n\tilde
a_{ij}v_{y_iy_j}+\sum_{i=1}\left( \tilde b_i-
\lambda\int_0^1\tilde\phi_i\,dz\right) v_{y_i}+h=0\\
v(y,T_0)=\tilde g(y),\end{array}\right.\end{equation} where
\[h(y,t)= \lambda(t)\int_0^1\left(v(y+\tilde\phi,t)-v(y,t)\right)\,dz\]
and $\tilde a_{ij}$ are the coefficients of the matrix
$\tilde\beta\tilde\beta^*/2$.

Proposition~3.3 in \cite{P} yields the estimate
\begin{equation}\label{est} \vert v(y,t)-v(\tilde y,\tilde
t)\vert\leq C\left((1+\vert y\vert)\vert t-\tilde t\vert^{1/2}+\vert
y-\tilde y\vert\right)\end{equation} for some constant $C$. Together
with the assumptions on $\phi$, this implies that $h\in
C_\alpha(\R^n\times [0,T_0])\cap C_{pol}(\R^n\times [0,T_0])$.
Consequently, Theorems~A.14 and A.18 in \cite{JT2} ensure the
existence of a unique classical solution $w\in C_{pol}(\R^n\times
[0,T_0])\cap C^{2,1}_\alpha(\R^n\times [0,T_0])$ to equation
(\ref{vis}). This classical solution $w$ to (\ref{vis}) can also be
represented (through the Feynman-Kac representation theorem) as
\[w(y,t)=E_{y,t}\left( \int_t^{T_0}h(Z(s),s)\,ds+\tilde g(Z(T_0))\right)\]
where $Z=(Z_1,...,Z_n)$ is the continuous diffusion process given by
\[dZ_i=\left( \tilde b_i(Z(t),t)-
\lambda(t)\int_0^1\tilde\phi_i(Z(t),t,z)\,dz\right)\,dt+
\sum_{j=1}^n\tilde\beta_{ij}(Z(t),t)\,dW_j\] and $Z(t)=y$. Since $h$
is Lipschitz continuous in $y$, it follows from Lemma~3.1 in
\cite{P} that $w$ is Lipschitz continuous in $y$, uniformly in $t$.
From the uniqueness result Theorem~4.1 in \cite{P} we deduce that
$v=w$. Consequently, $v\in C_{pol}(\R^n\times [0,T_0])\cap
C^{2,1}_\alpha(\R^n\times [0,T_0])$, and therefore $h\in
C^{2,0}_\alpha(\R^n\times[0,T_0])$. Applying Theorem~A.18 in
\cite{JT2} to equation (\ref{vis}) once again we find that $v=w\in
C^{4,1}_\alpha(\R^n\times [0,T_0])$. Transforming back to the
original coordinates we get $u\in C^{4,1}_\alpha(\R_+^n\times
[0,T_0])$.
\end{proof}

\section{The LCP condition as a necessary condition
for preservation of convexity}
\label{LCP condition}

Following \cite{JT2}, see also \cite{EJT} and \cite{ET}, to
investigate which models are convexity preserving we introduce the
notion of \emph{locally convexity preserving} (LCP) models.

\begin{definition}
An operator $\mathcal D$, where $\mathcal D$ equals either $\mathcal
M$, $\mathcal A$ or $\mathcal B$ specified above, is LCP at a point
$(x,t)\in\R_+^n\times [0,T]$ if for any direction
$v\in\R^n\setminus\{0\}$ we have that
\[\partial_v^2(\mathcal D f)(x,t)\geq 0\]
for all convex functions $f\in C^4_\alpha(\R_+^n)\cap
C^2_{pol}(\R_+^n)$ with $f_{vv}(x)=0$.
\end{definition}

We then have the following key result.

\begin{theorem}
\label{hekscher} Consider the following statements:
\begin{itemize} \item[(i)] The model is convexity preserving.
\item[(ii)] The operator $\mathcal M=\mathcal A+\mathcal B$ is LCP
at all points $(x,t)\in\R_+^n\times [0,T]$.
\item[(iii)]
The operators $\A$ and $\B$ are both LCP at all points
$(x,t)\in\R_+^n\times [0,T]$.
\end{itemize}
We have that (i)$\implies$(ii)$\iff$(iii).
\end{theorem}

\begin{remark}
Under a few additional growth conditions on $\beta$, $\phi$ and
their derivatives, it is possible to prove also (ii)$\implies$(i).
Thus all three statements in Theorem~\ref{hekscher} are equivalent,
see Theorem~4.3 in \cite{ET} for the one-dimensional case. We do not
pursue this further since we only use LCP as a necessary condition
in the analysis below.
\end{remark}

\begin{proof}
To prove (i)$\implies$(ii) we argue as in the proof of Lemma~3.3 in
\cite{JT2}. Choose $(x_0,T_0)\in\R_+^n\times[0,T]$ and let $g\in
C^4_\alpha (\R^n_+)\cap C^2_{pol} (\R^n_+)$ with $g_{vv}(x_0)=0$ for
some direction $v$. Define $u$ to be the solution to
\[u_t+\M u=0\]
on $\R^n_+\times [0,T_0]$ with terminal condition $u(x,T_0)=g(x)$.
Let $\tilde g\in C_\alpha^4(\R_+^n)$ be convex, Lipschitz continuous
and satisfy $\tilde g=g$ inside a box which contains $x_0$ and all
possible values of $x_0+\phi$. It follows that $\B g=\B \tilde g$ in
a (spatial) neighborhood of $x$, so $\partial_v^2 (\B
g)(x_0,T_0)=\partial_v^2(\B\tilde g)(x_0,T_0)$. Now, let $\tilde u$
be the solution to
\[\tilde u_t+\M \tilde u=0\]
on $\R^n_+\times [0,T_0]$ with terminal condition $u(x,T_0)=\tilde
g(x)$. Then, since the model is convexity preserving, $\tilde u$ is
convex in $x$ at all times prior to $T_0$, and in particular $\tilde
u_{vv}(x_0,t)\geq 0$. Consequently, $\partial_tu_{vv}(x_0,T_0)\leq
0$, so it follows from Theorem~\ref{reg} that \begin{eqnarray*} 0
&\leq& -\partial_t\tilde u_{vv}(x_0,T_0)=-\partial_v^2\tilde
u_t(x_0,T_0)=\partial_v^2\M \tilde u(x_0,T_0)\\ &=& \partial_v^2 (\M
\tilde g)(x_0,T_0)=
\partial_v^2 (\M g)(x_0,T_0),\end{eqnarray*}
which finishes the proof of (i)$\implies$(ii).

The implication (iii)$\implies$(ii) is immediate from the definition
of LCP. It remains to show that (ii)$\implies$(iii), i.e. that if
$\M$ is LCP, then both $\A$ an $\B$ are LCP. This follows from
Lemmas~\ref{lem1} and \ref{lem2} below. \end{proof}

\begin{lemma}\label{lem1}
Let $t_0\in[0,T]$ be given, and let $f\in C^4_\alpha( \R_+^n)\cap
C^2_{pol}( \R_+^n)$ be a convex function with $f_{vv}(x_0)=0$ at
some point $x_0$ and for some direction $v$. Then, for any $\ep>0$
there exists a convex function $h\in C^4_\alpha(\R_+^n)\cap
C^2_{pol}( \R_+^n)$ with $h_{vv}(x_0)=0$ such that $\partial_v^2(\A
h)(x_0,t_0)=0$ and $\left\vert\partial_v^2(\B
(f-h))(x_0,t_0)\right\vert\leq \ep$. Consequently, if $\M$ is LCP,
then also $\B$ is LCP.
\end{lemma}

\begin{proof}
Let $(x_0,t_0)\in\R_+^n\times [0,T]$ be given, and assume that $f\in
C^4_\alpha (\R_+^n)\cap C^2_{pol}( \R_+^n)$ is a convex function
with $f_{vv}(x_0)=0$ for some direction $v$. Without loss of
generality we assume that $f(x_0)=0$ and $\nabla f(x_0)=0$ (this can
be done since $\A \tilde f=\B \tilde f=0$ for all affine functions
$\tilde f$). Let $D$ be a constant such that for all unit directions
$w$ we have $f_{ww}(x)\leq D$ for all $x$ in a neighborhood of
$x_0$. Let $\varphi:\R_+^n\to[0,1]$ be a smooth function such that
\[\varphi(y)=\left\{\begin{array}{ll}0 & \mbox{ if }\vert y\vert\leq
1\\ 1 & \mbox{ if }\vert y\vert\geq 2,\end{array}\right.\] and let
$C_1$ be a constant such that
$\vert\varphi_w(y)\vert \leq C_1$ and $\vert
\varphi_{ww}(y)\vert\leq C_1$ for all $y\in\R_+^n$ and all unit directions
$w$. Further, let $\psi:[0,\infty)\to [0,\infty)$ be a smooth and
non-decreasing function satisfying
\[\psi(s)=\left\{\begin{array}{ll}
0 & \mbox{for }s\in[0,1/2]\\
Ms & \mbox{for }s\in[1,2]\\
3M & \mbox{for }s\geq 3,\end{array}\right.\] where $M$ is a constant
satisfying $M> 8C_1D$. Now, for $\delta>0$, define the function
$h=h^\delta$ by
\[h(x)=f\left(\varphi\left(\delta^{-1}(x-x_0)\right)
(x-x_0)+x_0\right)+\delta^2\int_0^{\vert x-x_0\vert/\delta}\psi(
s)\,ds.\] Then $h$ is the wanted function for some $\delta$ small
enough. Indeed, first note that $h$ is 0 if $\vert x-x_0\vert\leq
\delta/2$. Consequently, $\partial_v^2( \A h)(x_0,t_0)=0$. Moreover,
$h$ is convex if $\delta$ is small enough. To see that
$\partial_v^2(\B (f-h))$ can be made small, note that
\begin{eqnarray*}
\partial_v^2(\B f)(x_0,t_0)&=& \lambda\int_0^1\left(
\partial_v^2(f(x_0+\phi))-f_{vv}(x_0)-\partial^2_v(\phi\cdot\nabla
f(x_0))\right)\,dz\\
&=& \lambda\int_0^1\left(
\partial_v^2(f(x_0+\phi))-\phi_{vv}\cdot\nabla
f(x_0)\right)\,dz\end{eqnarray*} where we have used $f_{vv}(x_0)=0$
and $f_{vw}(x_0)=0$ for any direction $w$ (the latter statement
follows from $f_{vv}=0$ and the convexity of $f$). Similarly,
\[\partial_v^2(\B h)(x_0,t_0)
=\lambda\int_0^1\left(\partial_v^2(h(x_0+\phi))-\phi_{vv}\cdot
\nabla h(x_0)\right)\,dz.\] Thus, since $\nabla f(x_0)=\nabla
h(x_0)$,
\begin{eqnarray*}\partial_v^2(\B (f-h))(x_0,t_0)
&=&\lambda\int_0^1\partial_v^2\left(f(x_0+\phi)-h(x_0+\phi)\right)\,dz\\
&=& \lambda\int_0^1\partial_v^2 \left(f(x_0+\phi)-h(x_0+\phi)\right)
1_{\{z:\vert\phi\vert< 3\delta\}}\,dz\\ && +
\lambda\int_0^1\partial_v^2\left(\delta^2 \int_0^{\vert\phi\vert /
\delta}\psi(s)\,
ds\right)1_{\{z:\vert\phi\vert\geq 3\delta\}}\,dz\\
&=& I_1 + I_2.
\end{eqnarray*}
Here $I_1$ converges to 0 as $\delta$ goes to 0 since
$\phi(x,t,z)\not=0$ for almost all $z$ by (A6) and
$\partial_v^2(h(x+\phi))1_{\{z:\vert\phi\vert< 3\delta\}}$ is
bounded uniformly in $\delta$ at $x=x_0$. Similarly,
\begin{eqnarray*}
I_2 &=& \lambda\int_0^1\left(\delta\vert\phi\vert_{vv} \psi(\vert
\phi\vert/\delta) + \vert\phi\vert_v^2\psi^\prime (\vert
\phi\vert/\delta)\right)
1_{\{z:\vert\phi\vert\geq 3\delta\}}\,dz\\
&=& \lambda\int_0^1\delta\vert\phi\vert_{vv} 3M
1_{\{z:\vert\phi\vert\geq 3\delta\}}\,dz,
\end{eqnarray*}
so it follows from (A6) that also $I_2$ converges to 0. Thus, by
choosing $\delta$ small enough, we find that $\partial_v^2(\B
h)(x_0,t_0)$ is arbitrarily close to $\partial_v^2(\B f)(x_0,t_0)$.
\end{proof}

\begin{lemma}
\label{lem2} Let $t_0\in[0,T]$ be given, and let $f\in C^4_\alpha (
\R_+^n)\cap C^2_{pol}( \R_+^n)$ be a convex function with
$f_{vv}(x_0)=0$ at some point $x_0$ and for some direction $v$.
Then, for any $\ep>0$ there exists a convex function $h\in
C^4_\alpha(\R_+^n)\cap C^2_{pol}( \R_+^n)$ with $h_{vv}(x_0)=0$ such
that $\partial_v^2(\A h)(x_0,t_0)=\partial_v^2(\A f)(x_0,t_0)$ and
$\left\vert\partial_v^2(\B h)(x_0,t_0)\right\vert\leq \ep$.
Consequently, if $\M$ is LCP, then also $\A$ is LCP.
\end{lemma}

\begin{proof}
Without loss of generality, assume that $f(x_0)=0$ and $\nabla
f(x_0)=0$. As is seen in the proof of Theorem~\ref{wittbrattstrom}
below,
\begin{eqnarray}
\label{hermansson}
\partial_v^2(\B h)(x_0,t_0)&=&\int_0^1\Big(h_{vv}(x_0+\phi)+2\phi_v\cdot\nabla
h_v(x_0+\phi)\\
\notag
&&\hspace{10mm}+\phi_vHh(x_0+\phi)\phi_v^*\\
\notag && \hspace{10mm}+ \phi_{vv}\cdot\big(\nabla
h(x_0+\phi)-\nabla h(x_0)\big)\Big)\,dz
\end{eqnarray}
provided $h$ is convex and $h_{vv}(x_0)=0$ (here $Hh$ denotes the
Hessian of $h$). Thus it suffices to find $\delta>0$ and $h$
satisfying $h=f$ on $\{\vert x-x_0\vert\leq \delta\}$ and such that
$\nabla h$ and $H h$ are small on $\{\vert x-x_0\vert\geq \delta\}$.

To do this, let $C_1$ be a constant such that $f_{ww}\leq C_1$ in a
neighborhood of $x_0$ and for all unit directions $w$. Let
$\varphi:\R^n\to \R$ be a smooth non-negative function satisfying
$0\leq\varphi\leq 1$ such that
\[\varphi(x)=\left\{\begin{array}{cl}
1 & \mbox{if }\vert x-x_0\vert\leq 1\\
0 & \mbox{if }\vert x-x_0\vert\geq 2.\end{array}\right.\] Further,
let $C_2$ be a constant such that $\vert\varphi_w\vert\leq C_2$ and
$\vert\varphi_{ww}\vert\leq C_2$ for all unit directions $w$. Let
$M> 6C_1C_2$ and let $\psi:[0,\infty)\to[0,\infty)$ be a smooth and
non-decreasing function satisfying
\[\psi(s)=\left\{\begin{array}{cl}
0 & \mbox{for }s\in[0,1/2)\\
Ms & \mbox{for }s\in(1,2)\\
3M & \mbox{for }s\in(3,\infty).\end{array}\right.\] Now let
\[h^\delta (x):=f(x)\varphi((x-x_0)/\delta)+\delta^2\int_0^{\vert
x-x_0\vert/\delta}\psi(s)\,ds.\] Then $h=h^\delta$ is the wanted
function for some $\delta$ small enough. Indeed, first note that $\A
f(x_0,t_0)=\A h(x_0,t_0)$ since $h\equiv f$ for $\vert
x-x_0\vert\leq \delta$. Also note that straightforward calculations
yield that $h$ is convex if $\delta$ is small enough. Moreover
$h(x)=k\vert x-x_0\vert +b$ for $\vert x-x_0\vert\geq 3\delta$,
where $k=3M\delta$ and $b$ is some constant. If $w$ is a unit vector
with the same direction as $\phi$, i.e. if $\phi=\vert\phi\vert w$,
then it follows that
\[h_{x_i}(x_0+\phi)-h_{x_i}(x_0) = \int_0^{\vert\phi\vert}h_{x_iw}(x_0+sw)\,ds\]
can be made arbitrarily small (when varying $\delta$) since
$h_{x_iw}$ is bounded inside $\vert x-x_0\vert\leq 3\delta$
(uniformly in $\delta$) and vanishes outside this region. Thus the
last term in the right hand side of (\ref{hermansson}) can be made
arbitrarily small.

Moreover, examining the first three terms of (\ref{hermansson}) one
finds that these together form a second derivative of $h$, evaluated
at $x_0+\phi$, in the direction $w:=v+\phi_v$. Now
\[0\leq h_{ww}(x_0+\phi)\leq k\vert w\vert^2/\vert \phi\vert\]
for $\vert \phi\vert\geq 3\delta$. Since $k$ is linear in $\delta$
it follows that also the three first terms in (\ref{hermansson}) can
be made arbitrarily small when decreasing $\delta$.
\end{proof}

\section{A characterization of LCP models}
\label{LCP}

In \cite{EJT} it is shown that, within a large class of models, the
only differential operators of the form $\A$ which are LCP in
dimension $n\geq 2$ are the ones corresponding to geometric Brownian
motions. In that sense, there are not very many convexity preserving
diffusion models in higher dimensions. In this section we study the
LCP-condition for the operator $\B$ corresponding to the jump part
of $X$.

Let $Hf$ denote the Hessian of a function $f$. The following theorem
gives a precise description of which jump structures $\phi$ give
rise to an integro-differential operator $\B$ which is LCP.

\begin{theorem}
\label{wittbrattstrom} The operator $\B$ is LCP at a point $(x,t)$
if and only if for all directions $v\in\R^n\setminus\{ 0\}$ we have
\begin{eqnarray}
\label{ohly} &&\int_0^1\Big(f_{vv}(x+\phi)+2\phi_v\cdot\nabla
f_v(x+\phi)+\phi_vHf(x+\phi)\phi_v^*\\
\notag && \hspace{9mm} + \phi_{vv}\cdot\big(\nabla f(x+\phi)-\nabla
f(x)\big)\Big)\,dz\geq 0
\end{eqnarray}
for all convex functions $f\in C^4_\alpha(\R_+^n)\cap C^2_{pol}
(\R_+^n)$ with $f_{vv}(x)=0$.
\end{theorem}

\begin{proof}
Assume that $f\in C^4_\alpha(\R_+^n)\cap C^2_{pol}(\R_+^n)$ is
convex and that $f_{vv}=0$ at $x$. Straightforward calculations
yield that
\begin{eqnarray*}
\partial_v^2 (\mathcal B f)(x)
&=&
\int_0^1\Big(f_{vv}(x+\phi)+2\phi_v\cdot\nabla f_v(x+\phi)+\phi_v Hf(x+\phi)\phi_v^*\\
&& +\phi_{vv}\cdot\nabla f(x+\phi)-\phi_{vv}\cdot\nabla f(x) -
2\phi_v\cdot\nabla f_v -\phi\cdot\nabla f_{vv}\Big)\,dz.
\end{eqnarray*}
From the assumption $f_{vv}=0$ at $x$ it follows, due to the
convexity of $f$, that $\nabla f_v=\nabla f_{vv}=0$ at $x$. Thus
\begin{eqnarray*}
\partial_v^2 (\mathcal B f)(x)
&=&
\int_0^1\Big(f_{vv}(x+\phi)+2\phi_v\cdot\nabla f_v(x+\phi)+\phi_vHf(x+\phi)\phi_v^*\\
&& +\phi_{vv}\cdot\big(\nabla f(x+\phi)-\nabla f(x)\big)\Big)\,dz.
\end{eqnarray*}
Consequently, $\B$ is LCP if and only if (\ref{ohly}) holds for all
directions $v\in\R^n\setminus\{0\}$ and all $f\in C^4_\alpha
(\R_+^n)\cap C^2_{pol}(\R_+^n)$ with $f_{vv}(x)=0$.
\end{proof}

\begin{corollary}
\label{cor} Let $(x,t)\in\R_+^n\times [0,T]$. If \begin{equation}
\label{cond} \int_0^1 \phi_{vv}\cdot\big(\nabla f(x+\phi)-\nabla
f(x)\big)\,dz\geq 0
\end{equation}
for all convex functions $f\in C^4_\alpha (\R_+^n)\cap C^2_{pol}
(\R_+^n)$ and all directions $v\in \R^n\setminus\{ 0\}$, then the
operator $\B$ is LCP at $(x,t)$.
\end{corollary}

\begin{proof}
Since $f$ is convex, and since $f_{vv}+2\phi_v\cdot\nabla
f_v+\phi_vHf\phi_v^*$ is the second derivative of $f$ in the
$(v+\phi_v)$-direction, it is clear that (\ref{cond}) is sufficient
for the LCP-condition.
\end{proof}

\begin{remark}
If for all $i=1,...,n$, the function $\phi_i$ is convex in $x$ at
all points $(x,t,z)$ where $\phi_i(x,t,z)$ is positive, and $\phi_i$
is concave in $x$ at all points $(x,t,z)$ where $\phi_i(x,t,z)$ is
negative, then (\ref{cond}) is clearly satisfied. This sufficient
condition for preservation of convexity was used in \cite{ET} in a
one-dimensional setting. Also note that it is possible to show that
the condition (\ref{cond}) is strictly weaker than the condition
(\ref{ohly}).
\end{remark}

\section{The case of only finitely many possible jump sizes}
\label{finite} In this section we investigate models with only
finitely many possible jump sizes at each time. More specifically,
we assume that for each fixed $x$ and $t$, the function
$z\mapsto\phi(x,t,z)$ takes at most finitely many values.

\begin{theorem}
\label{werner} Assume there are only finitely many jump sizes, and
let $(x,t)\in\R_+^n\times [0,T]$. Then the following conditions are
equivalent:
\begin{itemize}
\item[(i)] $\B$ is LCP at $(x,t)$.
\item[(ii)]
(\ref{ohly}) holds for all directions $v$ and all convex $f\in
C^4_\alpha (\R_+^n)\cap C^2_{pol} (\R_+^n)$ with $f_{vv}=0$.
\item[(iii)]
(\ref{cond}) holds for all directions $v$ and all convex $f\in
C^4_\alpha (\R_+^n)\cap C^2_{pol} (\R_+^n)$.\end{itemize}
\end{theorem}

\begin{proof}
In view of Theorem~\ref{wittbrattstrom} and Corollary~\ref{cor} we
only need to show the implication (ii)$\implies$(iii).

To do this, let $f\in C^4_\alpha (\R_+^n)\cap C^2_{pol} (\R_+^n)$ be
convex. Since there are only finitely many possible values of
$\phi$, we can deform $f$ to be flat around all possible values of
$x+\phi$ and also in a neighborhood of $x$ without altering the
first derivative at these points. Accordingly, the first three terms
of the integrand in (\ref{ohly}) vanish, whereas the last term
remains unchanged, which finishes the proof.
\end{proof}

\begin{theorem}
\label{winberg} Assume there are only finitely many jump sizes and
that the model is convexity preserving. Then, for
each $i=1,...,n$, the ``jump volatility''
\[\sqrt{\frac{\lambda^2(t)\int_0^1\phi_i^2(x,t,z)\,dz}{x_i^2}}\]
of the $i$th asset is increasing as a function of $x_i$ at each fixed
time $t$.
\end{theorem}

\begin{proof}
Since the model is convexity preserving, it follows from
Theorem~\ref{hekscher} and Theorem~\ref{werner} that the inequality
(\ref{cond}) holds at all points and for all convex functions $f\in
C^4_\alpha(\R_+^n)\cap C^2_{pol}(\R_+^n)$. Choosing $f=x_i^2$ in
(\ref{cond}) gives that for all $i=1,...,n$ and for all directions
$v$ we have
\begin{equation}
\label{bildt} \int_0^1(\phi_i)_{vv}\phi_i\,dz\geq 0
\end{equation}
at all points $(x,t)$. Fix an $i$ and let
\[\psi(x,t):=\frac{\int_0^1\phi_i^2(x,t,z)\,dz}{x_i^2}.\]
Differentiating $\psi$ with respect to $x_i$ gives
\[x_i^4\psi_{x_i}=2\int_0^1\left(x_i^2\phi_i(\phi_i)_{x_i}-x_i\phi_i^2\right)\,dz.\]
Using (\ref{bildt}) with $v=e_{x_i}$ (here $e_{x_i}$ denotes
the $i$th unit coordinate
vector), integration by parts, $\phi_i\to 0$ as $x_i\to 0$ and
Jensen's inequality we find
\begin{eqnarray}\label{strindberg}
0 &\leq& x_i^2\int_0^{x_i}\int_0^1 (\phi_i)_{x_ix_i}\phi_i\,dzdx_i\\
&=& \notag\int_0^1 \left(x_i^2\phi_i(\phi_i)_{x_i}-
x_i^2\int_0^{x_i}(\phi_i)_{x_i}^2\, dx_i\right)\,dz\\ &\leq&
\notag\int_0^1 \left(x_i^2\phi_i(\phi_i)_{x_i}
-x_i\phi_i^2\right)\,dz.
\end{eqnarray}
This shows that $\psi_{x_i}$ is non-negative, which finishes the
proof.
\end{proof}

To the best of our knowledge, not very many models in finance have
increasing volatilities. Instead, models have typically large
volatilities for small values of the underlyings. If we restrict our
attention to these typical models, we show below that preservation
of convexity is a rather special property.

\begin{theorem}
\label{gBm} Let $n\geq 2$. Assume that there are only finitely many
possible jump sizes and that for all $i$ and $j$, $\beta_{ij}$ is a
function merely of $x_i$ and $t$. Also assume that for all
$i=1,...,n$ and for each fixed time $t$ the ``total volatility''
\[\sqrt{\frac{\beta_{i1}^2(x_i,t)+...+\beta_{in}^2(x_i,t)
+\lambda^2(t)\int_0^1\phi_i^2(x,t,z)\,dz}{x_i^2}}\] of the $i$th
asset is not an increasing function of $x_i$, unless it is constant.
If the model is convexity preserving, then $\beta_{ij}$ and $\phi_i$
are linear in $x_i$ for all $i$ and $j$, and $\phi_i$ does not
depend on $x_j$ for $j\not= i$. More explicitly, there exist
functions $\gamma_{ij}:[0,T]\to\R$ and $\gamma_i:[0,T]\times[0,1]$
such that
\begin{equation}\label{bet}
\beta_{ij}(x_i,t)=x_i\gamma_{ij}(t)\end{equation} and
\begin{equation}
\label{phi}\phi_i(x,t,z)=x_i\gamma_i(t,z)\end{equation} for almost
all $z$.
\end{theorem}

\begin{proof}
First note that according to Theorem~\ref{hekscher} both the
operators $\mathcal A$ and $\mathcal B$, corresponding to the
diffusion part and the jump part of $X$, respectively, are LCP. Now
note that if the total volatility is strictly decreasing in some
interval, then either a "diffusion volatility"
\[\sqrt{\frac{\beta_{i1}^2(x_i,t)+...+\beta_{in}^2(x_i,t)}{x_i^2}}\]
or a "jump volatility"
\[\sqrt{\frac{\lambda^2(t)\int_0^1\phi_i^2(x,t,z)\,dz}{x_i^2}}\]
is strictly decreasing in some interval. However, the proof of
Theorem~2.3 in \cite{EJT} implies that all diffusion volatilities
are increasing, and Theorem~\ref{winberg} above implies that all
jump volatilities are increasing. Consequently, all diffusion
volatilities and all jump volatilities are constant in $x_i$.

It then follows from the proof of Theorem~2.3 in \cite{EJT} that all
$\beta_{ij}$ are linear in $x_i$. Moreover, if the $i$th jump
volatility is constant in $x_i$, then the corresponding inequalities
in (\ref{strindberg}) reduce to equalities. Since Jensen's
inequality reduces to an equality if and only if the integrand is
constant, we find that for almost all $z$ the function $\phi_i$ has
to be linear in $x_i$. It thus only remains to show that $\phi_i$
does not depend on $x_j$, $j\not= i$. To do this we fix $j\not= i$,
and we plug $v=e_{x_i}+se_{x_j}$ and $\phi_i=\gamma_i(\overline
x,t,z)x_i$ into the inequality (\ref{bildt}), where $s\in\R$ and
$\overline x=(x_1,...,x_{i-1},x_{i+1},...,x_n)$. We find that
\[0 \leq \int_0^1\phi_i(\phi_i)_{vv}\,dz
= 2sx_i\int_0^1\gamma_i(\gamma_i)_{x_j}\,dz+s^2x_i^2\int_0^1
\gamma_i(\gamma_i)_{x_jx_j}\, dz.\] Since this expression is
non-negative for any choice of $s\in\R$, we must have
\begin{equation}\label{persson}\int_0^1\gamma_i(\gamma_i)_{x_j}\,dz=0.
\end{equation}
Now let $x_j^\prime\in\R_+$. Performing similar calculations as in
(\ref{strindberg}), but with $v=e_{x_j}$, we find that
\begin{eqnarray*}
0 &\leq& (x_j-x_j^\prime)^2\int_0^1\left(
\phi_i(x,t,z)(\phi_i)_{x_j}(x,t,z) -
\phi_i(x,t,z)(\phi_i)_{x_j}(x^\prime,t,z)\right)\,dz \\
&&  - \vert x_j-x_j^\prime\vert \int_0^1
\left(\phi_i(x,t,z)-\phi_i(x^\prime,t,z)\right)^2\,dz,
\end{eqnarray*}
where $x^\prime=(x_1,...,x_{j-1},x_j^\prime,x_{j+1},...,x_n)$. In
view of (\ref{persson}), the first integral vanishes. Consequently,
$\phi_i$ is $z$-almost surely constant in $x_j$, which finishes the
proof.
\end{proof}

\begin{remark}
Note that the models satisfying (\ref{bet}) and (\ref{phi}) are all
convexity preserving. Indeed, by explicit solution formulas,
\[X_i(T_0)=X_i(0)\exp\left\{-\frac{1}{2} \sum_{j=1}^n
\int_0^{T_0} \gamma_{ij}^2(t)\,dt
+\sum_{j=1}^n\int_0^{T_0}\gamma_{ij}(t)\,dW_j\right\}J_i(T_0),\]
where
\[J_i(T_0)=\exp\left\{-\int_0^{T_0}\int_0^1\lambda(t)\gamma_i(t,z)
\,dzdt\right\} \prod_{0\leq z\leq 1}\prod_{0\leq t\leq T_0}(1 +
\gamma_i(t))^{v(z,t)},\] so $X_i(T_0)$ is linear in the starting
value $X_i(0)$. Consequently,
\[u(x,0)=E_{x,0}g(X(T_0))\]
is convex in $X(0)$ provided $g$ is convex.
\end{remark}

\end{document}